\documentclass{article}
\usepackage{amssymb,amsmath,amsthm}
\newtheorem{thm}{Theorem}[section]
\begin{document}
\title{About one characteristic of set of linear
differential systems with non-negative coefficients}
\author{I.~Kopshaev
\thanks{Institute of mathematics of NAS of KAZAKHSTAN, 125 Pushkina str.,
050010 Almaty, KAZAKHSTAN email: {\em{kopshaev@math.kz}}} \and
A.~Sultanbekova \thanks{Institute of mathematics of NAS of
KAZAKHSTAN, 125 Pushkina str., 050010 Almaty, KAZAKHSTAN email:
{\em{sultanbekova@math.kz}}}}

\maketitle
\date
\begin{abstract}
The families of morphisms of vector fibre bundle (\cite{Mill1})
defined by the linear systems of differential equations with
non-negative coefficients is considered. Authors proved that the
specified families of morphisms is not saturated (\cite{Mill2}).
\end{abstract}

We consider the vector fibre bundle $(E,p,B)$ with $R^{n}$ as a
fibre and $B$ as a base (where $B$ is full metric space). On
$(E,p,B)$ has fixed some Riemannian metric (\cite{Husemoller}, P.
58-59).

Investigate the families of morphisms $\mathfrak{G}$ of linear
enlargement of dynamic system (\cite{Mill3}):
$$
(X(m), \chi(m)): (E,p,B) \to (E,p,B),
$$
($m \in N$), of the vector fibre bundle $(E,p,B)$, where
\begin{equation}
\label{kop_eq1} B=M_n, \quad E = B \times R^n, \quad p = pr_1,
\end{equation}
$$ X^t(A,x) = (\chi^tA, \mathfrak{X}(t,0,A) \cdot x),
$$
$$ \chi^t A(\cdot) = A(t+(\cdot)),$$
here $M_n$ - the space of linear systems of differential equations
with non-negative coefficients (\cite{Rahim}), $A \in B$, $x \in
R^n$, $\mathfrak{X}(\Theta, \tau, A)$ - Cauchy matrix of the
system $\dot{x} = A(t) \cdot x$.

\begin{thm}
The families of morphisms $\mathfrak{G}$ of the vector fibre
bundle (\ref{kop_eq1}) is not saturated.
\end{thm}

\bibliographystyle{plain}
\bibliography{nonexistence}

\end{document}